\documentclass{amsart}

\usepackage{amssymb}

\theoremstyle{plain}

\numberwithin{equation}{section}

\newtheorem{teo}[equation]{Theorem}

\newtheorem*{teorema}{Theorem A}
\newtheorem*{teoremab}{Theorem B}
\newtheorem*{teoremac}{Theorem C}

\newtheorem{cora}[equation]{Corollary}

\newtheorem{lem}[equation]{Lemma}

\newtheorem{prop}[equation]{Proposition}

\newcommand{\Irr}{\operatorname{Irr}}
\newcommand{\Lin}{\operatorname{Lin}}

\newcommand{\Ker}{\operatorname{Ker}}

\newcommand{\SL}{\operatorname{SL}}

\theoremstyle{definition}

\newtheorem{rem}[equation]{Remark}
\begin{document}	
\title{Squares of characters and finite groups}

\author{Edith Adan-Bante}

\address{University of Southern Mississippi Gulf Coast, 730 East Beach Boulevard,
 Long Beach MS 39560}

\email{Edith.Bante@usm.edu}

\keywords{Products of characters, groups of odd order, irreducible constituents}

\subjclass{20c15}

\date{2004}
\begin{abstract} Let $G$ be a group of odd order and $\chi$ be a complex irreducible
 character. Then there exists a unique character $\chi^{(2)}\in\Irr(G)$
such that $[\chi^2,\chi^{(2)}]$ is odd. Also, there exists a unique character 
$\psi\in \Irr(G)$
such that $[\psi^2, \chi]$ is odd. 
\end{abstract}
\maketitle

\begin{section}{Introduction}

Let $G$ be a finite group. Denote by $\Irr(G)$ the set of 
irreducible complex characters of $G$. 
Through this work,
we use the notation of \cite{isaacs}. In addition, we are going to 
denote by $\Lin(G)=\{\lambda\in \Irr(G) \mid \lambda(1)=1\}$ the set of
linear characters.

Let $\vartheta$ be a generalized character of $G$, i.e $\vartheta \in \mathbb{Z}[\Irr(G)]$.
We define $\vartheta^{(2)}$ by $\vartheta^{(2)}(g)= \vartheta(g^2)$ for all $g\in G$.
 By Theorem 4.5 of \cite{isaacs}, we have that
$\vartheta^{(2)}$ is also a  generalized character.  Also we can check that
 if $\vartheta\in \Irr(G)$ and
$G$ is of odd order, then $\vartheta^{(2)}\in \Irr(G)$ (see Exercise 4.5 of \cite{isaacs}). 

Let $G$ be a group of odd order and $\lambda\in \Lin(G)$ .
Since $\lambda^2\in \Lin(G)$,  $\lambda^2$ has a unique irreducible constituent 
that appears with odd multiplicity, namely $\lambda^2$,
 and such character has the same degree as
$\lambda$. Since $|G|$ is odd and $\lambda\in \Lin(G)$, we can check that 
there exists  $\mu\in \Lin(G)$ such that $\lambda$ is the unique 
irreducible constituent that appears with odd multiplicity in the character
$\mu^2$.
So we can find an irreducible character that is  ``the square" 
(``the square root")
 given any linear character. 
 For any character $\chi \in \Irr(G)$, where $\,\chi$ is 
 not necessarily linear,  the following result implies that 
$\chi^{(2)}\in \Irr(G)$ is the unique irreducible character of $G$
that appears with odd multiplicity in the character $\chi^2$. So we
can regard $\chi^{(2)}$ as ``the square" of $\chi$. 
	
\begin{teorema}

Let $G$ be a finite group of  odd order. 
The map $\vartheta \mapsto \vartheta^{(2)}$ is an automorphism of the 
ring $\mathbb{Z}[\Irr(G)]$. Also, if $\vartheta\in \mathbb{Z}[\Irr(G)]$
then 
\begin{equation}\label{square}
 \vartheta^2= \vartheta^{(2)} + 2 \Gamma
\end{equation}
\noindent  for some 
$\Gamma\in \mathbb{Z}[\Irr(G)]$.
\end{teorema}
  
  We now study squares of irreducible characters of prime degree.
 
 \begin{teoremab}\label{chi(1)=2}
 Let $G$ be a finite group and $\chi\in\Irr(G)$ be a monomial character such that 
 $\chi(1)=2$.  Then one of the following holds
  
(i) $\chi^2=\lambda_1 +\lambda_2 +\theta$, for some $\lambda_1, \lambda_2\in \Lin(G)$
 and some $\theta\in \Irr(G)$ with $\theta(1)=2$.
 
(ii) $\chi^2=\sum_{i=1}^4\lambda_i$ for some  distinct linear
characters $\lambda_1, \lambda_2, \lambda_3$ and $\lambda_4$.
\end{teoremab}

 Let $\Theta$ be a character of $G$. Denote by
 $\eta(\Theta)$ the number of distinct irreducible constituents of
 the  character $\Theta$. In Theorem \ref{morethan1} is proved that if 
$G$ is  a supersolvable group and $\chi\in \Irr(G)$ with $\chi(1)=2^n$ for some integer 
$n>0$, then $\eta(\chi^2)\geq 2$, that is 
$\chi^2$ is the sum of at least two distinct irreducible
 characters. 
 \begin{teoremac}
Let $p$ be an odd prime number.
 Let $G$ be a finite nilpotent group and $\chi\in\Irr(G)$  
 with $\chi(1)=p$. Then all the irreducible constituents of $\chi^2$ have degree 
 $p$ and  one of the following holds:

(i) $\eta(\chi^2)=1$. Also $\chi^2=p\chi^{(2)}$ and 
 $\chi$ vanishes outside ${\bf Z}(\chi)$.

(ii) $\eta(\chi^2)=\frac{p+1}{2}$. Furthermore
$$\chi^2=\chi^{(2)} + 2\sum_{i=1}^{\frac{p-1}{2}}\theta_i$$
\noindent  where  $\theta_1, \theta_2, \ldots, \theta_{\frac{p-1}{2}}$ are 
distinct irreducible characters of $G$.
 \end{teoremac}
 
  Let $G$ be a nilpotent group and $\chi\in \Irr(G)$ with $\chi(1)=p$.
 By Theorem B and Theorem C we have that  
 if $p$ is a prime number, then there are exactly 2
 possible values for $\eta(\chi^2)$.
 Assume now that $\chi(1)=p^n$ for a fixed prime $p$. 
 The author
 wonders  for what integer $n>1$, the number of possible values for $\eta(\chi^2)$
 really does depend upon the prime $p$.

In Remark \ref{super}, examples are provided, showing that Theorem C
may not remain true if $G$ is only a  
supersolvable group of odd order.

{\bf Acknowledgment.} I thank Professor Everett C. Dade for 
very useful conversations and 
emails.
I also thank Professor Avinoam Mann and the referee for their corrections and 
suggestions to improve this note. 

\end{section}

\begin{section}{Proofs}
  
\begin{proof}[Proof of Theorem A]
Exercise 4.7 of \cite{isaacs} implies \eqref{square}. It remains then to prove that
the map $\vartheta \mapsto \vartheta^{(2)}$ is a bijection.

Since $G$ has odd order,  
 given any 
 character $f$ and any element $g\in G$, there exists a unique element $f'(g)$ over the complex 
such that $(f'(g))^2=f(g)$. 
Let $f'$ be the function defined as
$(f'(g))^2=f(g)$ for all
$g\in G$. If $f$ is any character, so is $f'$ and $(f')^2=f$. 
 That shows that the map $f\mapsto f'$ is onto on the set of characters,
  and therefore one to one.
 We can then
extend this to generalized characters. 
\end{proof}
   
Let $G$ be the dihedral group of order $8$ and 
$\chi\in \Irr(G)$ with $\chi(1)=2$.
We can check that 
$$\chi^2=\sum_{\lambda \in \Lin(G)} \lambda.$$  
Thus for a group $G$ of even order and $\chi\in \Irr(G)$, there may not exist a unique
character that appears with odd multiplicity in $\chi^2$. 
\begin{proof}[Proof of Theorem B]
Without lost of generality, we may assume that $\chi$ is a faithful character.
Since $\chi$ is monomial, there exists a subgroup $N$ of $G$
 and a linear character $\lambda$ of $N$ such that $\lambda^G=\chi$.
Since $\chi(1)=2$, then $|G:N|=2$ and thus $N$ is a normal subgroup of $G$.
  Since
$\chi$ is faithful and $\chi(1)=2$, we have that $N$ is an abelian subgroup of $G$. 
Fix $g\in G\setminus N$.
By Clifford theory we have that $\lambda^g\neq\lambda$ and  
$\chi_N=\lambda+\lambda^g$.
Thus
$$(\chi^2)_N=\lambda^2+(\lambda^g)^2+2\lambda\lambda^g.$$
Observe that $\lambda\lambda^g$ is a $G$-invariant character of $N$. Since $|G:N|=2$,
it follows that $\lambda\lambda^g$ extends to $G$. 

Let $\mu\in\Lin(G)$ be
an extension of $\lambda\lambda^g$. Observe that 
$\chi\nu=\chi$ for any $\nu\in \Lin(G/ N)$. Thus $\chi^2$ has at least 2 
linear constituents, namely the two  distinct extensions of
the character $\lambda\lambda^g$. By Exercise 4.12 of \cite{isaacs}, we have that
$[\chi^2,\theta]\leq \theta(1)$ for any $\theta\in \Irr(G)$. Thus $\chi^2$ has at least 3 distinct irreducible constituents, i.e. $\eta(\chi^2)\geq 3$.
Since $\chi^2(1)=4$, either
$\eta(\chi^2)=3$ or $\eta(\chi^2)=4$. 
 
 If $\eta(\chi^2)=3$ then $\chi^2$ has 2 linear constituents and one irreducible
 constituent of degree 2 and (i) holds.
 
 If $\eta(\chi^2)=4$ then all the irreducible constituents of $\chi^2$ are linear
 and (ii) holds.
\end{proof}

\begin{rem} Let $\SL(2,3)$ be the special linear group of degree 2 over the field of 3 elements, and $\chi\in \Irr(\SL(2,3))$ with 
$\chi(1)=2$. We can check that $\chi^2=\lambda + \theta$, where $\theta\in \Irr(G)$ and
$\theta(1)=3$ and $\lambda\in \Lin(\SL(2,3))$.
 Thus Theorem B does not remain true without the hypothesis that 
$\chi$ is a monomial character.
\end{rem}
\begin{rem} 
Let $C_{2^n}$ and $C_2$ be cyclic subgroups of order $2^n$, where $n>1$, and $2$ respectively.
Set $N=C_{2^n}\times C_{2^n}$. Observe that $C_2$ acts on $N$ by permutation. 
Let $G$ be the semidirect product of $N$ and $C_2$. Thus $|G|=2^{2n+1}$.

Let $(\lambda,  1_{C_{2^n}}   )\in \Irr(N)$ be a character such that
$\lambda^2\neq 1_{C_{2^n}}$. Set $\chi= (\lambda,  1_{C_{2^n}})^G$ and 
$\theta= (\lambda^2,  1_{C_{2^n}})^G$.
Observe that $\chi$ and $\theta$ are irreducible characters of degree 2.
  We can check that $$\chi^2=\lambda_1+\lambda_2+\theta$$
\noindent where
$\lambda_1$ and $\lambda_2$ are the 2 distinct extensions of the $G$-invariant 
character
$(\lambda,\lambda)\in \Lin(N)$.

Let $G$ be the dihedral subgroup of order 8. We can check that $\chi^2$ is the sum
of 4 linear characters.

We conclude that for each of the cases in Theorem B, there exists 
a group 
$G$ and a character $\chi\in \Irr(G)$ satisfying the conditions in that case. 
\end{rem}

\begin{teo}\label{morethan1}
Let $G$ be a finite supersolvable group of $G$ and $\chi\in \Irr(G)$. Assume that
$\chi(1)=2^n$ for some integer $n>0$. Then $\eta(\chi^2)>1$. 
\end{teo}
\begin{proof}
Assume that the statement is false. Let  $G$ be a supersolvable group and  
$\chi\in \Irr(G)$ be a character of $G$ such that they are a  minimal 
counterexample with respect to the order of $G$. So $\chi^2=m\theta$ 
for some integer $m$, some $\theta\in \Irr(G)$ and $\chi(1)=2^n$ for some 
integer $n>0$. 
   
Let $Z={\bf Z}(\chi)$ be the center of the character $\chi$. Since $\chi^2=m\theta$,
we can check that $Z={\bf Z}(\theta)$. Let $\lambda\in\Irr(Z)$ be the unique 
linear character of $Z$ such that $\chi_Z=\chi(1) \lambda$. Observe that 
$\theta_Z=\theta(1)\lambda^2$ since $\chi^2=m\theta$. 
 Let $Y/Z$ be a chief factor of $G$ and $\mu\in \Irr(Y)$ be a character such that 
 $[\chi_Y, \mu]\neq 0$. Because $G$ is supersolvable and $Y/Z$ is a chief factor 
 of $G$, it follows that $Y/Z$ is cyclic and therefore $\lambda$ extends to $Y$. Thus 
 $\mu$ is a linear character.
 
  Let $H=G_{\mu}$ be the stabilizer of $\mu$ in $G$ and $\chi_{\mu}\in \Irr(H)$
   be the Clifford
  correspondent of $\chi$ and $\mu$, i.e $\chi_{\mu}^G=\chi$. 
 Since $\mu$ is a linear character, $[\chi_Y,\mu]\neq 0$  and $Y>Z={\bf Z}(\chi)$,
 we have that $H$ is a proper subgroup of $G$. Observe that $[\theta_Y, \mu^2]\neq 0$. 
 We can check that the stabilizer $G_{\mu^2}$ of $\mu^2$ in $G$ contains
 $G_{\mu}=H$ and $G_{\mu^2}<G$. Let $K$ be a subgroup 
 of $G$ such that $ G_{\mu^2}\leq K<G$ and  $\chi_{\mu}^K$ is not a linear character. 
 Observe such subgroup $K$ exists since $\chi(1)>2$ by Theorem B.
 
 By Clifford theory we have that $\chi_{\mu}^K\in \Irr(G)$. Since $\chi_{\mu}^K$ lies
 above $\mu$, 
 all the irreducible constituents of $\chi_{\mu}^K\chi_{\mu}^K$ lie above $\mu^2$. 
 Since $ G_{\mu^2}\leq K$, by Clifford theory we can check
  that $\eta( \chi_{\mu}^K\chi_{\mu}^K)=
 \eta((\chi_{\mu}^K\chi_{\mu}^K)^G)$. Since $|K|<|G|$ and $\chi_{\mu}(1)$ is not a linear
 character,   we have 
 that $\eta( \chi_{\mu}^K\chi_{\mu}^K))>1$. Since $(\chi_{\mu}^K)^G=\chi$, we can check that 
 $\eta((\chi_{\mu}^K\chi_{\mu}^K)^G)\leq \eta(\chi\chi)$.
 But then $1< \eta( \chi_{\mu}^K\chi_{\mu}^K)\leq \eta(\chi\chi)=1$. Therefore
 such group $G$ and such character $\chi$ can not exist. 
 \end{proof}	
 Recall that the irreducible characters of supersolvable 
 groups are monomial. 
The author wonders if Theorem \ref{morethan1}
 remains true with the weaker hypothesis that 
$\chi$ is monomial.
  
 A corollary of Theorem \ref{morethan1} is the following

\begin{cora} 
Let $G$ be a nilpotent group and $\chi\in \Irr(G)$. If $\chi^2$ is a multiple of 
an irreducible, then $\chi(1)$ is an odd integer. 
\end{cora}

 \begin{prop}\label{lemma1}
Let $G$ be a group of odd order  and $\chi\in \Irr(G)$. 
  If  $\chi(1)>1$  then $[\chi^2,\lambda]=0$ for all $\lambda \in \Lin(G)$. 
\end{prop}
\begin{proof}
 Assume that $[\chi^2,\lambda]\neq 0$ for some $\lambda\in \Lin(G)$.
Since $G$ has odd order, 
there exists some character $\beta\in \Lin(G)$ such that $\beta^2=\lambda$. 
Thus $[\chi^2,\lambda]=[\chi^2, \beta^2]=[\chi\overline{\beta}, \overline{\chi}\beta]$.
Since $\chi\overline{\beta}, \overline{\chi}\beta \in \Irr(G)$, it follows that 
$$\chi\overline{\beta}= \overline{\chi}\beta.$$
Observe that
$\overline{\chi}\beta=\overline{\chi\overline{\beta}}$. 
 Since $\chi\overline{\beta}\in \Irr(G)$ is a real character 
and $G$ has odd order,
 it follows that  $\chi\overline{\beta}=1_G$. Thus $\chi=\beta$ and 
$\chi(1)=1$. 
 
 \end{proof}
\begin{prop}\label{lemma2}
Let $G$ be a finite group  and $\chi\in \Irr(G)$. 
Assume that $\Ker(\chi)\neq {\bf Z}(\chi)$.  Then $[\chi^2,\chi]=0$. 
\end{prop}

\begin{proof}
Observe that $[\chi^2,\chi]=[\chi,\chi\overline{\chi}]$. Observe also
that $\Ker(\chi \overline{\chi})={\bf Z}(\chi)$. Thus if $[\chi,\chi\overline{\chi}]\neq 0$, then $\Ker(\chi)\geq {\bf Z}(\chi)$. Therefore $\Ker(\chi)={\bf Z}(\chi)$.
\end{proof}

\begin{cora}
Let $G$ be a finite $p$-group, where $p$ is a prime number, and $\chi\in \Irr(G)$.
If $\chi\neq 1_G$ then $[\chi^2,\chi]=0$.
\end{cora}
 
 \begin{rem} Let $D$ be the dihedral group of order $8$ and $\chi\in \Irr(D)$ be such that
 $\chi(1)=2$. We can check that $\chi^2$ is the sum of 4 linear characters. 
 Thus Proposition \ref{lemma1}  may not hold true when $G$ has even order. 

 Let $G$ be the semi-direct product $ <\tau> \rtimes <\sigma>$, where
$\sigma$ and $\tau$ have orders $7$ and $3$, respectively, and
$\sigma^{\tau} = \sigma^2$. We can check that any faithful irreducible
character $\chi$ of $G$ has degree $3$, ${\bf Z}(\chi)=1$ and 
$[\chi^2, \chi] = 1$.
Thus Proposition \ref{lemma2} may not remain true if $\Ker(\chi)={\bf Z}(\chi)$.
 \end{rem}
\begin{lem}\label{chi(1)=p}
 Let $G$ be a finite $p$-group, where $p$ is an odd prime, and 
$\chi\in \Irr(G)$ with $\chi(1)=p$. Then all the irreducible constituents of $\chi^2$ have degree $p$ and  one of the following holds:

(i) $\eta(\chi^2)=1$. Also $\chi^2=p\chi^{(2)}$ and 
 $\chi$ vanishes outside ${\bf Z}(\chi)$.

(ii) $\eta(\chi^2)=\frac{p+1}{2}$. Furthermore
$$\chi^2=\chi^{(2)} + 2\sum_{i=1}^{\frac{p-1}{2}}\theta_i$$
\noindent  where  $\theta_1, \theta_2, \ldots, \theta_{\frac{p-1}{2}}$ are 
distinct irreducible characters of $G$. 
\end{lem}

\begin{proof}
 By working with the group $G/\Ker(\chi)$, without lost of generality 
we may assume that $\Ker(\chi)=1$. 

Assume that $\chi^2$ is a multiple of an irreducible.
In \cite{edithma} it is proved that if $\chi, \psi\in \Irr(G)$ are faithful characters
of a finite $p$-group $G$, and  the product
$\chi\psi $ is a multiple of an irreducible,  then  $\chi$ and $\psi$ vanish
outside the center of $G$. Since $\chi^{(2)}$ is a constituent of $\chi^2$ and $\chi^2$ is a 
multiple of an irreducible, (i)
follows.

We may assume now that $\chi^2$ is not a multiple of an irreducible. 
In \cite{edith} it is proved that if $\chi,\psi\in \Irr(G)$ are faithful characters
of a finite $p$-group, where $p$ is odd, either $\chi\psi$ is a multiple of an
irreducible or $\chi\psi$ has at least $\frac{p+1}{2}$ distinct irreducible
constituents. By Proposition \ref{lemma1} we have that $\chi^2$ does not have
any linear constituent. Since $\chi^2(1)=p^2$ and $\chi^2$ is not a multiple of 
an irreducible, it follows that all the constituents of $\chi^2$ have degree
$p$. By Theorem A we have that $\chi^2$ has a unique irreducible constituent
that appears with odd multiplicity.
Since $\chi^2$ has at least $\frac{p+1}{2}$ distinct irreducible
constituents, $\chi^{(2)}$ is the 
unique irreducible constituent that appears with odd multiplicity and
$$\chi(1)^2=p^2= p+ 2(\frac{p-1}{2})p,$$
\noindent 
it follows
that $\chi^2$ has to have exactly $\frac{p+1}{2}$ distinct irreducible constituents
and so (ii) follows.

 \end{proof}
\begin{proof}[Proof of Theorem C]
Since any nilpotent subgroup is isomorphic to the direct sum of its $p$-Sylow subgroups,
 Theorem C follows from Lemma \ref{chi(1)=p}.
 \end{proof}
\begin{rem}\label{super}
Let $p$ and $q$ odd prime numbers such that $q=rp+1$, for some integer
$r$, with $1<r<\frac{p+1}{2}$.
Let $G$ be a nonabelian group of order $pq$. We can check that $G$ has exactly 
$p$ distinct linear characters and $r$ distinct irreducible characters of
degree $p$.  

 Let $\chi\in \Irr(G)$ with $\chi(1)=p$. Since $G$ has odd order, then 
$\chi^2$ does not have any linear constituent. Thus $\chi^2$ is the sum of
 at most $r$ distinct irreducible characters of $G$, where $r<\frac{p+1}{2}$.
We can check that  $\chi^2$ is not a multiple of an irreducible. Therefore
Theorem C does not remain true with the weaker hypothesis that $G$ is 
supersolvable of odd order. 
\end{rem} 

\end{section}

\end{document}